\documentclass[onefignum,onetabnum]{siamart171218}

\usepackage{amsmath,amsfonts,amssymb, color, ifpdf}
\usepackage{hyperref}

\graphicspath{ {./figs/} } 

\renewcommand{\vec}[1]{\boldsymbol{#1}}

\usepackage{enumerate}
\usepackage{subcaption}
\usepackage{cite}
\usepackage{booktabs}

\usepackage{tikz}          % figures
\usepackage{tikzscale}     % scaling figures with includegraphics
\usetikzlibrary{shapes,arrows,positioning,decorations.pathreplacing,
	intersections, matrix,patterns,fadings,calc,decorations.markings,arrows.meta}
\usepackage{pgfplots}

\newcommand{\Div}[1]{\nabla\cdot #1}
\newcommand{\Grad}[1]{\nabla #1}
\newcommand{\Curl}[1]{\nabla\times #1}

\renewcommand{\vec}[1]{\mathbf{#1}}
\renewcommand{\Re}{Re}
\newcommand{\Ha}{H\!a}%Without the \!, the kerning is very poor...

\newcommand{\beps}{\boldsymbol{\varepsilon}}

\newcommand{\dX}{\delta\vec{x}}
\newcommand{\dU}{\delta\vec{u}}
\newcommand{\dB}{\delta\vec{B}}
\newcommand{\dP}{\delta p}
\newcommand{\dR}{\delta r}

\newcommand{\intO}{\int_{\Omega}}

\newcommand{\intdx}{\;\mathrm{d}\vec{x}}

\DeclareMathOperator{\CurlOp}{curl}

\newcommand{\Hcurl}{H(\text{curl})}

\newsiamremark{remark}{Remark}
\newsiamremark{weakform}{Weak Form}

\headers{Monolithic~MG~for~MHD}{J.~H.~Adler, T.~Benson, E.~C.~Cyr, P.~E.~Farrell, S.~MacLachlan, R.~Tuminaro}
\title{Monolithic Multigrid for
    Magnetohydrodynamics\thanks{Submitted to the editors
      DATE.\funding{This research was enabled in part by support provided by ACENET (\url{www.ace-net.ca}), Scinet (\url{www.scinethpc.ca}), and Compute Canada (\url{www.computecanada.ca}).  The work of JA was partially supported by National Science Foundation grants DMS-1216972 and DMS-1620063.  The work of TB was performed under the auspices of the U.S. Department of Energy by Lawrence Livermore National Laboratory under Contract DE-AC52-07NA27344 (LLNL-JRNL-811901).  The work of
        EC and RT was supported by the U.S.~Department of Energy, Office of Science, Office of Advanced Scientific Computing Research, Applied Mathematics program.  Sandia National Laboratories is a multimission laboratory managed and operated by National Technology and Engineering Solutions of Sandia, LLC., a wholly owned subsidiary of Honeywell International, Inc., for the U.S. Department of Energy's National Nuclear Security Administration under grant DE-NA-0003525.  The work of PEF was supported by the UK Engineering and Physical Sciences
      Research Council [grant number EP/R029423/1] and by a London Mathematical Society `Research in Pairs' grant. The work of SM was
        partially supported by an NSERC Discovery Grant.  This paper describes objective technical results and analysis. Any subjective views or opinions that might be expressed in the paper do not necessarily represent the views of the U.S. Department of Energy or the United States Government.}}}
  \author{James H.~Adler\thanks{Department of Mathematics, Tufts University, 503 Boston Ave., Medford, MA 02155 (\email{james.adler@tufts.edu}).}
    \and Thomas Benson\thanks{Lawrence Livermore National Laboratory, 7000 East Avenue, P.O. Box 808, L-495, Livermore, CA 94551 (\email{benson31@llnl.gov})}
    \and Eric C.~Cyr\thanks{Center for Computing Research, Sandia National Labs, Albuquerque, NM 87185-1318 (\email{eccyr@sandia.gov}), Livermore, CA 94551 (\email{rstumin@sandia.gov})}
    \and Patrick E.~Farrell\thanks{Mathematical Institute, University of Oxford, Oxford, UK (\email{patrick.farrell@maths.ox.ac.uk}).}
    \and Scott MacLachlan\thanks{Department of Mathematics and
      Statistics, Memorial University of Newfoundland, St.\ John's,
      NL, Canada (\email{smaclachlan@mun.ca}).}
    \and Ray Tuminaro\footnotemark[4]}
\date{\today}

\begin{document}

\maketitle

\begin{abstract}
  The magnetohydrodynamics (MHD) equations model a wide range of
  plasma physics applications and are characterized by a nonlinear
  system of partial differential equations that strongly couples a
  charged fluid with the evolution of electromagnetic fields. After
  discretization and linearization, the
  resulting system of equations is generally difficult to solve due to
  the coupling between variables, and
  the heterogeneous coefficients induced by the linearization process. In this
  paper, we investigate multigrid preconditioners for this system
  based on specialized relaxation schemes that properly address the
  system structure and coupling.  Three extensions of Vanka relaxation
  are proposed and applied to problems with up to 170 million degrees of
  freedom and fluid and magnetic Reynolds numbers up to 400 for stationary problems and up to 20,000 for time-dependent problems.
\end{abstract}

\begin{keywords}
Magnetohydrodynamics (MHD), monolithic multigrid, Vanka relaxation
\end{keywords}

\begin{AMS}
 65F10, 65N55, 65N22, 76W05
\end{AMS}

\section{Introduction}

%\textcolor{red}{
%Standard spiel(?): block-factorization based preconditioners have
%become prominent for solving multiphysics problems, with lots of
%allusions to ``physics-based'' preconditioning, but monolithic MG has,
%when developed fully, generally proven to give better solvers.  Our
%goal: study this for the 4-field formulation of single fluid,
%viscoresistive MHD.
%}

The magnetohydrodynamics (MHD) equations model the evolution of a
quasi-neutral plasma in the presence of magnetic fields. This model is highly
nonlinear, and contains strong coupling between the fluid unknowns and the
magnetic field. As a result, the system is characterized by a number of varying
temporal and spatial scales. The temporal scales give rise to very
stiff modes that make explicit time integration methods intractable due
to stability limitations. To overcome this, unconditionally stable implicit time integration methods
are
often used \cite{2002ChaconL_KnollD_FinnJ-aa,BPhilip_LChacon_MPernice_2008a,JNShadid_EtAl_2010a,EGPhillips_EtAl_2014a,MaHuHuXu_2016a,EGPhillips_etal_2016a,shadid2016scalable,MWathen_CGreif_2020a}. These time integration approaches lead to large
sparse linear systems whose solution requires effective parallel preconditioners.

This
paper develops new preconditioners for linear systems arising from Newton linearization of a finite element discretization
of the viscoresistive,
incompressible MHD model. Previous work
\cite{JHAdler_TRBenson_EtAl_2016a} utilized a vector-potential
formulation to enforce the solenoidal constraint. Motivated by three-dimensional calculations, 
we use a formulation that maintains the
primitive variables, including the magnetic field, $\vec{B}$, and
enforces the solenoidal constraint weakly by using a Lagrange
multiplier that is added to Faraday's law
\cite{ASchneebeli_DSchotzau_2003a,DSchotzau_2004a,CGreif_EtAl_2010a}. Note
that the Lagrange multiplier approach is only one such technique for
ensuring the solenoidal condition is satisfied; others include exact-penalty
methods \cite{EGPhillips_EtAl_2014a}, using augmented Lagrangian terms \cite{HuQiuShi_2020a}, vector-potential formulations
\cite{JHAdler_TRBenson_EtAl_2016a,ECCyr_EtAl_2013a,JNShadid_EtAl_2010a,
JAdler_etal_2019b}, and compatible discretizations \cite{MaHuHuXu_2016a}.  The main contribution of this paper is a new family of Vanka-style relaxation schemes for a multigrid procedure that is appropriate for this type of Lagrange multiplier enforcement.

In order to solve the system
numerically, we use a mixed finite element discretization and Newton's method, resulting in
linearized problems of saddle-point type. Saddle-point problems arise in
various contexts, ranging from economics to fluid dynamics; an
extensive general overview can be found in
\cite{MBenzi_GHGolub_JLiesen_2005a}. Numerous solution and
preconditioning strategies have been proposed, and often the solver or
preconditioner is very closely tied to the problem being
solved.
% However, there are two common methods of preconditioning:
% block-factorization approaches
% \cite{HElman_EtAl_2008a,HCElman_EtAl_2005a,MurRehman_EtAl_2011a,RVerfurth_1984a}
% and monolithic multigrid
% \cite{SPVanka_1986a,SPMacLachlan_CWOosterlee_2011a,JHAdler_TRBenson_EtAl_Sub2015b,JHAdler_TRBenson_EtAl_2016a}.
Block-factorization approaches manipulate the constituent blocks of
the Jacobian operator in order to resolve the coupling in the system
using inner solvers
on auxiliary problems
\cite{murphy2000,HElman_EtAl_2008a,MurRehman_EtAl_2011a,RVerfurth_1984a}.
These techniques have been extensively studied and, by utilizing this
paradigm, solvers for simpler PDEs
can be exploited to create scalable preconditioners for
more complex problems, such as MHD
\cite{ECCyr_EtAl_2013a,MaHuHuXu_2016a}.  In this community, such approaches are often referred to as ``physics-based'' preconditioners, which construct approximations of the Schur complements
based on the timescales and particular stiff modes that the simulation must resolve
\cite{LChacon_DAKnoll_2006a,BPhilip_LChacon_MPernice_2008a}.

%\textcolor{red}{
%Refs to past work (ours and others) on monolithic MG, developing out
%of the fluids context.  Focus on Vanka and Braess-Sarazin as competing
%options for relaxation (and Uzawa as a ``simplified Braess-Sarazin''
%that some folks like for reasons that remain unclear to me), coming
%from the fluids end.  Discuss relation between these and curl-curl
%relaxation schemes (Hiptmair and AFW), hinting there's more to see
%here than has been exposed.
%}

An alternative approach is to apply multigrid monolithically to the coupled linear system, exploiting the
structure of the PDE in the relaxation method. To this end, we propose a monolithic
multigrid approach for the MHD discretization presented in \cite{DSchotzau_2004a}.
Monolithic multigrid for coupled
systems is, in fact, one of the earliest ideas in the multigrid
literature \cite{ABrandt_1984a,ABrandt_NDinar_1979a}. Monolithic
schemes for incompressible fluid flow problems are well-studied, with
several treatments of the Stokes equations
\cite{CWOosterlee_FJGaspar_2008a,FGaspar_EtAl_2014a,JSchoberl_WZulehner_2003a}
and of the Navier-Stokes equations
\cite{MLarin_AReusken_2008a,SPVanka_1986a,VJohn_GMatthies_2001a,VJohn_LTobiska_2000a}. In
the MHD literature, a monolithic nonlinear multigrid method is investigated
for a finite-difference discretization in
\cite{MFAdams_RSamtaney_ABrandt_2010a}. A fully-coupled AMG approach
for a vector-potential formulation of resistive MHD is shown in
\cite{JNShadid_EtAl_2010a}, and a Lagrange multiplier MHD formulation using nodal
finite elements in~\cite{shadid2016scalable}.  However, this latter approach 
relies upon an
equal-order discretization in which unknowns for each variable are
collocated at mesh vertices; thus, it cannot be used for the mixed
discretization shown here. Similarly, fully-coupled AMG methods have
been used in the context of first-order system least squares discretizations of MHD
\cite{JAdler_EtAl_2010b}. Finally, in
\cite{JHAdler_TRBenson_EtAl_2016a}, two families of relaxation methods
for the vector-potential formulation are described and used within a
monolithic geometric multigrid preconditioner. In this paper, variations on the Vanka-style
relaxation methods are extended to the Lagrange multiplier
formulation described below.

Following the formulation of the MHD equations in
\cite{ASchneebeli_DSchotzau_2003a,DSchotzau_2004a}, we introduce a
Lagrange multiplier, $r$, to weakly enforce the solenoidal constraint,
$\Div\vec{B}=0$. This Lagrange multiplier appears as a nonphysical
term, $\Grad r$, in Faraday's law. The set of viscoresistive MHD
equations in a bounded Lipschitz domain $\Omega\subset\mathbb{R}^d$ ($d \in \{2, 3\}$) that we consider in
this paper is then written,
\begin{align}
\label{eqn:ns}\partial_t {\bf u} -\Div\left[\frac{2}{\Re}\beps(\vec{u})\right] +
  (\vec{u}\cdot\Grad)\vec{u} + \Grad p - (\Curl\vec{B})\times\vec{B}
  &= \vec{f}, \\
\label{eqn:faraday}\partial_t {\bf B} + \frac{1}{\Re_m}\Curl\Curl\vec{B} -
  \Curl(\vec{u}\times\vec{B}) - \Grad r &= \vec{g}, \\
\label{eqn:incomp}\nabla\cdot{\bf u} & = 0, \\
\label{eqn:div-b-free}\nabla\cdot{\bf B} & = 0, 
\end{align}
where the unknown fields are the velocity ${\bf u}$, the magnetic
field ${\bf B}$, the pressure $p$ and the Lagrange multiplier $r$.  The strain-rate tensor is
$\beps(\vec{u}) = \frac{1}{2}\left(\Grad\vec{u} +
  \Grad\vec{u}^{\top}\right)$.
The nondimensional parameters are $\Re$, the hydrodynamic Reynolds
number, and $\Re_m$, the magnetic Reynolds number. The system is
closed with some appropriate set of boundary conditions dependent on the problem.  Preconditioners for this formulation were introduced in \cite{2017M.-WathenC_SchotzauD-aa, MWathen_CGreif_2020a, LLi_WZheng_2017a}.  In \cite{2017M.-WathenC_SchotzauD-aa, LLi_WZheng_2017a}, block-triangular preconditioners for the system are proposed and analyzed, while \cite{MWathen_CGreif_2020a} proposes a block-structured approximate inverse preconditioner.  In contrast, the monolithic multigrid approach considered here avoids the need to approximate one or more 
Schur complements (which may be a non-trivial task
for complex MHD systems) 
as required
%
%direct approximations of the Schur complements 
within block-factorization approaches.
%For concreteness, we consider Dirichlet boundary conditions on the velocity, tangential Dirichlet conditions on the magnetic field, and homogeneous Dirichlet conditions on the Lagrange multiplier in the exposition below.

% Note that the
% Lagrange multiplier approach
% (also used in \cite{CGreif_EtAl_2010a,PHouston_DSchotzau_XWei_2009a,wgs_2017})
% is only one such technique for including the solenoidal condition;
% others include exact-penalty methods \cite{EGPhillips_EtAl_2014a},
% vector-potenital formulations
% \cite{JNShadid_EtAl_2010a,ECCyr_EtAl_2013a,JHAdler_TRBenson_EtAl_Sub2015},
% and others.

\section{Discretization and linearization}

%\textcolor{red}{Ref to Sch\"otzau paper on mixed FEM discretization of MHD, with the
%``extra'' Lagrange multiplier (and Greif-Sch\"otzau-Michael Wathen
%work on block preconditioning for it), contrast with the vector
%potential formulation that we've used in the past, as well as other
%formulations in recent work from Eric et al.
%Probably some discussion of linearize-then-discretize and
%discretize-then-linearize, at least to make it clear that both are
%possible and are used, but we make a fixed choice (of the former,
%right?)}

System \eqref{eqn:ns}-\eqref{eqn:div-b-free} is
nonlinear, and we thus employ Newton's method as a nonlinear
solver.  We use a discretize-then-linearize approach and discretize the system using a mixed finite element method.  For simplicity, we first consider the steady-state version of
\eqref{eqn:ns}-\eqref{eqn:div-b-free}.
In the time-dependent case, we use a
second-order backward difference formula for the time derivative, and
then proceed with the discretization and linearization
described below. Additionally, we consider homogeneous Dirichlet boundary conditions on the velocity, homogeneous tangential Dirichlet conditions on the magnetic field, and homogeneous Dirichlet conditions on the Lagrange multiplier in this section's exposition.  To define the variational form, we consider the
appropriate solution spaces for each unknown:  
\begin{align*}
\vec{V} &= \vec{H}^{1}_0(\Omega) =
\{\vec{v}\in\vec{H}^{1}(\Omega)\;\colon\, \vec{v} = 0
\text{ on }\partial\Omega\},\\
\vec{C} &= \vec{H}_0(\CurlOp,\Omega) = \{\vec{c}\in \vec{L}^{2}(\Omega)\;\colon\,
\Curl{\vec{c}}\in\vec{L}^{2}(\Omega),\ \vec{n}\times\vec{c} = 0 \text{ on }\partial\Omega\},\\
Q &= L^{2}_0(\Omega) = \{q\in L^{2}(\Omega)\;\colon\, \intO q \intdx = 0\},\\
S &= H^{1}_0(\Omega) = \{s\in H^{1}(\Omega)\;\colon\, s = 0 \text{ on }\partial\Omega\},
\end{align*}
where $\vec{n}$ is the unit outward normal on the boundary.
We denote by $\|\cdot\|_0$ the usual norm on $L^{2}(\Omega)$ or the
vector version $\vec{L}^{2}(\Omega)$. Similarly, we denote by
$\|\cdot\|_1$ the usual norm in $H^{1}(\Omega)$ or
$\vec{H}^{1}(\Omega)$. Finally, for a vector
$\vec{c}\in\vec{H}(\CurlOp,\Omega)$, we define
$\|\vec{c}\|_{\CurlOp}^2 = \|\vec{c}\|_0^2 + \|\Curl\vec{c}\|_0^2$ to
be the norm in $\vec{H}(\CurlOp,\Omega)$.

Multiplying the MHD equations by the appropriate test functions drawn from the spaces above, 
then integrating over the domain and applying integration
by parts yields the \emph{nonlinear} variational form: Find
$(\vec{u},\vec{B},p,r)\in \vec{V}\times\vec{C}\times Q\times S$ such
that
\begin{align}
\intO \left (\frac{2}{\Re}\beps(\vec{u}):\beps(\vec{v}) - p\Div\vec{v}
+ \left[(\vec{u}\cdot\Grad)\vec{u}\right]\cdot\vec{v} -  \left[(\Curl\vec{B})\times\vec{B}\right]\cdot\vec{v}\right )\intdx\qquad ~\label{eqn:lin_var_momentum}\\
 = \intO \vec{f}\cdot\vec{v} \intdx
 % + \int_{\partial\Omega_N} p_{N}\vec{n}\cdot\vec{v}\intdS
 , \notag\\
\intO \left[\frac{1}{\Re_m}(\Curl\vec{B}) - (\vec{u}\times\vec{B})\right]\cdot(\Curl\vec{c}) -
\Grad r\cdot\vec{c}\intdx = \intO \vec{g}\cdot\vec{c}\intdx,\label{eqn:lin_var_faraday}\\
\intO q\Div\vec{u}\intdx = 0, \label{eqn:lin_var_continuity}\\
\intO \Grad s \cdot \vec{B}\intdx = 0,\label{eqn:lin_var_solenoid}
\end{align}
for all $(\vec{v},\vec{c},q,s)\in\vec{V}\times\vec{C}\times Q\times S$.
The boundary terms in \eqref{eqn:lin_var_faraday} vanish because
we strongly enforce that
$\vec{c}\in\vec{C} = \vec{H}_0(\CurlOp,\Omega)$, so that
$\vec{n}\times\vec{c} = 0$. Likewise, the boundary integral in
\eqref{eqn:lin_var_solenoid} has been eliminated by enforcing
$s\in S = H^1_0(\Omega)$.

We discretize
\eqref{eqn:lin_var_momentum}--\eqref{eqn:lin_var_solenoid} by approximating the solution of the
system using finite element functions
$(\vec{u}^{h},\vec{B}^{h},p^{h},r^{h})\in\vec{V}^{h}\times\vec{C}^{h}\times
Q^{h}\times S^{h}$,
where $\vec{V}^{h}\times Q^{h}$ is a standard inf-sup stable pair for
the incompressible Navier-Stokes problem, $\vec{C}^{h}$ is the first
family of N\'ed\'elec elements \cite{JCNedelec_1980a}, and $S^{h}$ is
an $H^1$-conforming space
\cite{ASchneebeli_DSchotzau_2003a,DSchotzau_2004a}. The domain $\Omega$ is 
partitioned into simplicial elements, and
thus we choose $\vec{P}_{2}-P_{1}$ (Taylor-Hood) elements for
$\vec{V}^{h}\times Q^{h}$. We choose the lowest-order N\'ed\'elec
space for $\vec{C}^{h}$ and $P_{1}$ for $S^{h}$. Well-posedness of both
the continuous and discrete formulations (under a suitable small-data assumption) is shown in
\cite{ASchneebeli_DSchotzau_2003a,DSchotzau_2004a}.
%, and block-triangular
%preconditioning techniques for the linearized version of this system were proposed in \cite{2017M.-WathenC_SchotzauD-aa}.

Linearizations of \eqref{eqn:lin_var_momentum}--\eqref{eqn:lin_var_solenoid} are computed
via Newton's method. Given an initial guess, $\vec{x}_0$, at the $n^{\text{th}}$ step,
we solve
\[
J(\vec{x}_n) \delta\vec{x} = - R(\vec{x}_n),
\]
where $\vec{x}_n = (\vec{u}^h_n,\vec{B}^h_n,p^h_n,r^h_n)^{\top}$ is the current
approximation to the solution,
$\dX = \vec{x}_{n+1} - \vec{x}_n = (\dU^h, \dB^h, \dP^h, \dR^h)^{\top}$ is the Newton
update,
$J(\vec{x}_n)$ is the Jacobian operator evaluated at $\vec{x}_n$, and
$R(\vec{x}_n)$ is the nonlinear residual evaluated at
$\vec{x}_n$.  Note that we assume the initial guess, $\vec{x}_0 = (\vec{u}^h_0,\vec{B}^h_0,p^h_0,r^h_0)^{\top}$, satisfies the appropriate boundary conditions for each unknown.

%\begin{remark}In this paper, we consider linearizing the system first,
%and then discretizing the resulting linearized system via the
%finite element method.  We refer to this as a
%\emph{linearize-then-discretize} approach.  Alternatively, in a
%\emph{discretize-then-linearize} approach, we would discretize the
%nonlinear equations directly and then use a nonlinear solver, such as
%the Full Approximation Storage (FAS) algorithm \cite{1984BrandtA-aa,2000BriggsW_HensonV_McCormickS-aa} or a Newton-Krylov-Schwarz method \cite{1998CaiX_GroppW_KeyesD_MelvinR_YoungD-aa}.  While the latter scheme can be more robust, the former is easier to use in practice and both methods are equivalent in the limit where the iterations approach the true solution \cite{2003CoddA_ManteuffelT_McCormickS-aa}.  
%Additionally, linearize-then-discretize approaches are very common in solving nonlinear PDEs.  Newton-Krylov methods, for instance, have been developed for a wide array of PDE systems and tuned for many applications (e.g.  \cite{2002ChaconL_KnollD_FinnJ-aa,2008PhilipB_ChaconL_PerniceM-aa,2010ShadidJ_PawlowskiR_BanksJ_ChaconL_LinP_TuminaroR-aa,1994EisenstatS_WalkerH-aa,1997ShadidJ_TuminaroR_H.-F.-WalkerH-aa,1998PerniceM_WalkerH-aa,2006PawlowskiR_ShadidJ_SimonisJ_H.-F.-WalkerH-aa}). 
%\end{remark}

Thus, the
sequence of discrete, linear variational problems arising from Newton's method leads to a
linear system of the following block form for each Newton step:
\begin{equation}\label{eqn:linear_system}
\mathcal{A}x =
\begin{bmatrix}
F & Z & B^{\top} & 0  \\
Y & D & 0   & C^{\top}\\
B & 0 & 0   & 0  \\
0 & C & 0   & 0
\end{bmatrix}\begin{bmatrix}
x_{\vec{u}} \\
x_{\vec{B}} \\
x_{p}      \\
x_{r}
\end{bmatrix} = \begin{bmatrix}
f_{\vec{u}} \\
f_{\vec{B}} \\
f_{p}      \\
f_{r}
\end{bmatrix} = b.
\end{equation}
Here, $x_{\vec{u}}$, $x_{\vec{B}}$, $x_{p}$, $x_{r}$ are the
discretized Newton corrections for $\vec{u}$, $\vec{B}$, $p$, and $r$,
respectively, and $f_{\vec{u}}$, $f_{\vec{B}}$, $f_{p}$, and $f_{r}$
are the corresponding blocks of the nonlinear residual.

\section{Monolithic MG}

Since the matrix in Equation \eqref{eqn:linear_system} is not symmetric, we will use
FGMRES \cite{YSaad_MHSchultz_1986a, YSaad_2003a} as the outer Krylov
method in a Newton-Krylov-multigrid approach to solving \eqref{eqn:linear_system}. Noting that \eqref{eqn:linear_system} is a saddle-point
problem, we develop an effective monolithic multigrid
preconditioner that naturally treats the structure of the block linear
system, avoiding the construction of approximate Schur complements
for the constraint degrees of freedom that are typically required in scalable block preconditioning approaches. As in any multigrid
method, the choice of a suitable relaxation scheme and its
complementarity with the coarse-grid correction is critically important to
achieve scalable performance. In this paper, we fix a robust
geometric coarse-grid correction procedure, and extend the well-known
Vanka relaxation scheme\cite{VJohn_LTobiska_2000a, SPVanka_1986b,
  MLarin_AReusken_2008a} to \eqref{eqn:linear_system}.  Specifically,
we generate a hierarchy of meshes by uniform refinement of a coarse
triangulation of the domain, generating $2^d$ fine-grid
elements from one coarse-grid element in $d$ dimensions. The interpolation operators
are chosen to be block-structured operators with diagonal blocks given
by the standard finite element
interpolation operators for each
variable,
\[
P = \begin{bmatrix}
P_{\vec{u}} \\
& P_{\vec{B}} \\
& & P_p \\
& & & P_r
\end{bmatrix},
\]
where $P_{\vec{u}}$ is the vector-quadratic ($\vec{P}_2$)
interpolation operator; $P_{\vec{B}}$ is the lowest-order first-family
N\'ed\'elec interpolation operator; and $P_p$ and $P_r$ are both the
linear ($P_1$) interpolation operator. Coarse-grid operators are
constructed by rediscretization.

\subsection{Vanka Relaxation}\label{ssec:vanka}
As is common when considering monolithic multigrid for a
fluid-dynamics problem, here we consider Vanka-type relaxation methods
\cite{SPVanka_1986a}. In \cite{JHAdler_TRBenson_EtAl_2016a}, these
methods were extended to a finite element discretization of a
vector-potential formulation of the MHD system, where only a single
constraint variable (the fluid pressure) is present. Here, we aim to
further extend Vanka relaxation schemes to the case of \eqref{eqn:linear_system}. A key difference here is the
presence of two Lagrange multipliers, as well as use of the full
$\vec{B}$ field, which is now discretized using curl-conforming vector
elements.

We begin with a review of the ideas of classical Vanka relaxation when
applied to the (Navier-)Stokes equations.  In that setting, we have a
natural partitioning of the sets of degrees of freedom (DoFs) in the
discrete operator into those associated with the (vector) velocity
approximation, $\mathcal{S}_{\vec{u}}=\{u_1,\dots,u_{n_u}\}$, and
those associated with the (scalar) pressure approximation,
$\mathcal{S}_{p}=\{p_1,\dots,p_{n_p}\}$.  Let $\mathcal{S} =
\mathcal{S}_{\vec{u}} \cup \mathcal{S}_{p}$ be the set of all the DoFs
in the system.  A ``Vanka relaxation scheme'' is typically taken to be
a block overlapping Gauss-Seidel (multiplicative Schwarz) iteration,
where the set of DoFs is partitioned, $\mathcal{S} = \cup_{\ell=1}^N
\mathcal{S}_{\ell}$, into $N$ ``Vanka blocks'' that are defined based
on the structure of the mesh and coupled system to be solved.
In fluid dynamics applications, the standard approach to decomposing
$\mathcal{S}$ into the subsets, $\mathcal{S}_{\ell}$, is to ``seed''
the choice of the Vanka blocks by the incompressibility constraint, or
(equivalently) by the pressure degrees of freedom.

Algebraically, this means that given a block linear system of the form, 
\begin{align*}
\mathcal{A} &=\begin{bmatrix} F & B^{\top} \\ B &
0 \end{bmatrix},
\end{align*}
we associate one Vanka block
to each row in the ``constraint'' matrix, $B$, and we define
$\mathcal{S}_{\ell}$ to be the set of degrees of freedom corresponding
to (symbolic) nonzero entries in row $\ell$ of $B$, as well as the seed DoF,
$p_{\ell}\in\mathcal{S}_{p}$.  Topologically, this can be seen as
isolating the nodal pressure DoFs (located at vertices of the mesh since
we use a $P_1$ discretization) and including all of the velocity DoFs
in the stencil surrounding them (in the closure of the star of the vertex). While
Vanka originally proposed such a block construction for the
Marker-and-Cell (MAC) staggered finite-difference discretization
scheme \cite{SPVanka_1986b}, its extension to finite element
discretizations has been considered in \cite{VJohn_LTobiska_2000a,
  SManservisi_2006a, MLarin_AReusken_2008a, HWobker_STurek_2009a,
  VJohn_GMatthies_2001a}. In these works, there is general agreement
that it is critical to form the Vanka blocks in order to preserve the
structure of the constraint blocks in the system. It is this
philosophy that we adopt below.

Once the Vanka blocks are formed, standard Vanka relaxation proceeds
by iterating (multiplicatively) over the blocks, sequentially solving
subproblems obtained by restricting the global linear system to the
block DoFs.
That is, for each Vanka block
$\mathcal{S}_{\ell}$, the global solution is updated according to
\begin{equation}\label{eqn:vanka_update}
x \leftarrow x + V_{\ell} \left( \omega {M}_{\ell\ell}^{-1}\right)V_{\ell}^{\top} \left(b -
\mathcal{A}x\right),
\end{equation}
and the updates are computed in a Gauss-Seidel fashion. Here,
$V_{\ell}^{\top}$ is a restriction operator that takes global vectors to
local vectors containing only the entries corresponding to DoFs in
block $\ell$, and $V_{\ell}$ is a prolongation operator that takes the
entries in a local vector over $\mathcal{S}_\ell$ and inserts them
appropriately into a global vector \cite{YSaad_2003a}. Typically, an
underrelaxation parameter, $\omega$, is used to improve performance of
the overall multigrid scheme.  Here, we take $\omega = 1$ within
the iteration and use an outer Chebyshev polynomial to accelerate the
relaxation scheme globally. Finally, $M_{\ell\ell}$ is the Vanka
submatrix (of dimension $|\mathcal{S_\ell}|\times|\mathcal{S_\ell}|$)
that is used to compute the update.  A common choice is to take
$M_{\ell\ell} = V_\ell^{\top}\mathcal{A}V_\ell$, the restriction of the
global system matrix to the DoFs in $\mathcal{S}_\ell$.  In the
finite-difference setting, substantial cost savings can often be realized by further economizing at this stage, by using
sparser approximations to the restricted global problem.  So-called
``diagonal Vanka'', for example, replaces the $F$ block in
$\mathcal{A}$ by its diagonal, $\text{diag}(F)$, defining
\[
M_{\ell\ell} = V_\ell^{\top}\begin{bmatrix} \text{diag}(F) & B^{\top} \\ B &
0 \end{bmatrix} V_\ell.
\]
Since a matrix of this form can be factored at low cost, such schemes
have been considered in several settings in the literature
\cite{JHAdler_TRBenson_EtAl_2016a, CRodrigo_etal_2016a, VJohn_LTobiska_2000a, SPMacLachlan_CWOosterlee_2011a}.
However, timing tests in \cite{JHAdler_TRBenson_EtAl_2016a} find there
is less advantage to this approach in the finite element setting,
since the cost of relaxation is largely dominated by assembling the
local residuals (computing $V_\ell^{\top}(b-\mathcal{A}x)$); for this
reason, we do not consider these approaches in detail in this work.

\subsection{Arnold-Falk-Winther relaxation for N\'ed\'elec elements}

Just as Vanka relaxation is well-established for problems in
computational fluid dynamics, a standard relaxation approach for
edge-element discretizations in $\Hcurl$ is that proposed
by Arnold, Falk, and Winther \cite{DNArnold_EtAl_2000a}. This relaxation (and the associated multigrid
method) is typically considered for finite element
discretizations of the standard second-order operators in these
spaces, such as finding $\vec{v}$ such that
\[
\nabla\times \left(\alpha \nabla\times \vec{v}\right) +
  \beta\vec{v} = \vec{g},
\]
for given parameters (possibly functions) $\alpha > 0$ and $\beta > 0$.  Arnold, Falk, and Winther proposed an overlapping
Schwarz relaxation scheme for such problems, collectively relaxing all DoFs in suitable patches
defined around each vertex, specifically the star of each vertex.  For regular 2D meshes and the lowest-order N\'{e}d\'{e}lec elements, as we consider here, this amounts to collecting the 6 edge DoFs connected to each node in each patch.

In our setting, ignoring the nonlinearity, we apply a variation on the Arnold-Falk-Winther (AFW) idea to the discretized form of \eqref{eqn:faraday}.  Here, however, we have a different lower-order term,
$\nabla\times(\vec{u}\times\vec{B})$, and the equation is modified by
the presence of a Lagrange multiplier, $r$.  Neither, however, requires any
fundamental changes in the structure of the relaxation. The Lagrange multiplier is incorporated by adding a single constraint DoF into each of the blocks prescribed in the AFW relaxation scheme (as the $P_1$ Lagrange multiplier is discretized with nodal finite elements, we add the DoF at the node central to the edges in each AFW patch).  A slight variation to the AFW
relaxation would be to view it in the same algebraic setting as
discussed for Vanka above, by taking each row of the discretized
constraint matrix, $C$ (from \eqref{eqn:linear_system}), and forming the relaxation patch by taking all
$\vec{B}$ DoFs that appear in that row as well as the Lagrange
multiplier DoF associated with it.  This results in the patch consisting of the closure of the star around each vertex, including all edge DoFs associated with elements that contain the constraint seed DoF, as illustrated in the center diagram of Figure \ref{fig:Vanka_blocks}.  In what follows, we use this larger patch in all variants proposed.

\subsection{Additive and multiplicative updates}

	Historically, Vanka iteration has generally been applied in a multiplicative
	manner, but it is also possible to apply it additively. Additive approaches
	are more straightforwardly parallelizable, at the cost of some degradation in convergence.
	Another significant advantage of additive methods is that they avoid the
	repeated computation of local residuals that is necessary in a multiplicative
	scheme; this was found to be the dominant cost of relaxation in \cite{JHAdler_TRBenson_EtAl_2016a}.
	For this reason, we consider only additive relaxation methods in Section \ref{sec:numerical_results}.

        Many variants of additive overlapping domain decomposition schemes exist in the literature, including those where overlap is accounted for by weighting corrections using a partition of unity and restricted additive Schwarz (RAS) methods.  Here, we consider only unweighted additive Schwarz methods, where the corrections in \eqref{eqn:vanka_update} are simply summed over each Vanka block, giving
\[
x \leftarrow x + \omega \sum_{\ell}V_{\ell} \left( \omega {M}_{\ell\ell}^{-1}\right)V_{\ell}^{\top} \left(b -
\mathcal{A}x\right),
\]
as the corresponding stationary iteration.  While other variants of additive Schwarz may lead to faster convergence of the domain-decomposition method as a standalone preconditioner, we emphasize that we use it only as a relaxation scheme in a monolithic multigrid method and find that this variant is sufficient for that task.

\subsection{Extending Vanka to the four-field formulation}\label{sec:vankaextensions}

In our view, there are three natural extensions of the Vanka
relaxation scheme described above that can be applied to the
discretized and linearized MHD system in \eqref{eqn:linear_system}.

\subsubsection{Segregated Vanka} The easiest approach to applying the above ideas to
\eqref{eqn:linear_system} is to do so in a {\it segregated} manner,
where we separately apply Vanka relaxation to the fluid-dynamics
(velocity-pressure) subproblem and AFW relaxation (extended, as discussed above, to the closure of the star of the vertex, including the Lagrange multiplier) to the
electromagnetics (magnetic field-Lagrange multiplier) subproblem.
This approach ignores the coupling in the $Y$ and $Z$ blocks of the
matrix in \eqref{eqn:linear_system} and, as such, is expected to lead
to degradation in performance for larger values of the hydrodynamic
and magnetic Reynolds numbers, where the coupling in these blocks is
important to resolving the solution of the discretized system.  We
refer to this relaxation scheme as {\it segregated Vanka} in the
numerical results below.  The grouping of degrees of freedom in this relaxation scheme is shown in Figure \ref{fig:Vanka_blocks}, where the classical Vanka coupling for the fluid degrees of freedom is shown on the left and the (extended) Arnold-Falk-Winther relaxation scheme for the electromagnetics subproblem is shown in the center.  One iteration of the segregated Vanka relaxation scheme is composed of independent additive relaxation schemes over both of these subproblems.

\subsubsection{Purist Vanka} In line with the classical fluid-dynamics relaxation, a second
approach is to consider a block two-by-two form of the matrix \eqref{eqn:linear_system},
\begin{equation}\label{eqn:sub_linear_system}
\mathcal{A}x = \begin{bmatrix}
\hat{F} & \hat{B}^{\top} \\
\hat{B} & 0
\end{bmatrix} \begin{bmatrix}
x_{\hat{\vec{u}}} \\
x_{\hat{p}}
\end{bmatrix} = \begin{bmatrix}
f_{\hat{\vec{u}}} \\
f_{\hat{p}}
\end{bmatrix},
\end{equation}
where we have grouped the vectors as
\[
x_{\hat{\vec{u}}} = \begin{bmatrix}
x_{\vec{u}} \\ x_{\vec{B}}
\end{bmatrix}, \quad
x_{\hat{p}} = \begin{bmatrix}
x_p \\ x_r
\end{bmatrix}, \quad
f_{\hat{\vec{u}}} = \begin{bmatrix}
f_{\vec{u}} \\ f_{\vec{B}}
\end{bmatrix}, \quad
f_{\hat{p}} = \begin{bmatrix}
f_p \\ f_r
\end{bmatrix},
\]
and the matrices as
\[
\hat{F} = \begin{bmatrix}
F & Z \\
Y & D
\end{bmatrix}, \quad
\hat{B} = \begin{bmatrix}
B & 0 \\ 0 & C \end{bmatrix}.
\]
In {\it purist Vanka}, we treat the block two-by-two constraint
matrix, $\hat{B}$, as a single constraint matrix, and define Vanka
blocks algebraically by looping row-by-row over $\hat{B}$, forming
Vanka blocks by taking all nonzero entries in each row of this matrix
along with the pressure or Lagrange-multiplier DoF associated with the
row.

The first problem with this approach is that, in our context, it leads to a relaxation scheme 
identical to the segregated Vanka procedure. This is 
due to the PDE structure as there are no nonzero entries associated
with magnetic field DoFs in rows of $\hat{B}$ associated with a
pressure DoF, and no velocity DoFs in the rows associated with a Lagrange-multiplier DoF.  
To avoid this limitation, instead of gathering the DoFs using the sparsity pattern of the
constraint matrix, we 
use the topological adjacency relation.  Since there is a nontrivial overlap
between each of the nodal pressure basis functions with magnetic field
basis functions associated with edges around it, and between the
nodal Lagrange multiplier basis functions and those of velocity, this
leads to greater coupling in the blocks.  Unfortunately, this coupling
leads to an additional challenge, since the submatrices associated
with the Vanka blocks constructed in this way are singular unless a
Lagrange multiplier DoF has been included in the patch.  
For the blocks containing only a pressure constraint, the singularity arises because
the discrete curl-curl operator has local
nullspace components, mirroring the fact that applying the curl to the gradient
of any scalar function is identically zero.
We address
this by using an augmented Lagrangian approach, which has been proposed
and used in \cite{GreifSchotzau2007}.  For Vanka blocks associated
with pressure degrees-of-freedom, we form the submatrices by first
adding a mass matrix on the N\'ed\'elec space to the $D$ block in the system matrix, then restricting to the
Vanka block.
The degrees of freedom in each resulting Vanka block match those depicted on the right of Figure \ref{fig:Vanka_blocks}, taking one of the two central nodal DoFs in each block.  Thus, the total number of blocks in this approach is equal to twice the number of vertices in the mesh.

\subsubsection{Coupled Vanka} 
\begin{figure}
	\centering
	\begin{tikzpicture}[scale=0.85]
	%\draw[help lines] (0,0) grid (13,4);
	
	%% Fluid Patch
	\draw[very thick] (0,0) -- (4,0) -- (4, 4) -- (0, 4) -- (0,0);
	\draw[very thick] (0,2) -- (4,2);
	\draw[very thick] (2,0) -- (2,4);
	\draw[very thick] (4,0) -- (0,4);
	\draw[very thick] (2,0) -- (0,2);
	\draw[very thick] (4,2) -- (2,4);
	\foreach \x in {1,...,3}{
		\foreach \y in {1,...,3}{
			\node[circle,draw,fill=gray!40,inner sep=0pt,minimum size=12pt] at (\x,\y) {}; 
		}
	}
	\node[circle,draw,fill=gray!40,inner sep=0pt,minimum size=12pt] at (2,0) {};
	\node[circle,draw,fill=gray!40,inner sep=0pt,minimum size=12pt] at (2,4) {};
	\node[circle,draw,fill=gray!40,inner sep=0pt,minimum size=12pt] at (0,2) {};
	\node[circle,draw,fill=gray!40,inner sep=0pt,minimum size=12pt] at (4,2) {};
	\node[circle,draw,fill=gray!40,inner sep=0pt,minimum size=12pt] at (4,0) {};
	\node[circle,draw,fill=gray!40,inner sep=0pt,minimum size=12pt] at (3,0) {};
	\node[circle,draw,fill=gray!40,inner sep=0pt,minimum size=12pt] at (4,1) {};
	\node[circle,draw,fill=gray!40,inner sep=0pt,minimum size=12pt] at (0,3) {};
	\node[circle,draw,fill=gray!40,inner sep=0pt,minimum size=12pt] at (0,4) {};
	\node[circle,draw,fill=gray!40,inner sep=0pt,minimum size=12pt] at (1,4) {};
	\node[circle,fill=black,inner sep=0pt,minimum size=7pt] at (2,2) {}; 
	\node at (1.6,1.6) {$p$};
	
	%% E&M Patch
	\tikzset{->-/.style={decoration={markings,mark=at position .6 with {\arrow[scale=1.7]{stealth}}},postaction={decorate}}}
	\draw[very thick] (5,0) -- (7,0);
	\draw[very thick,->-] (7,0) -- (9,0);
	\draw[very thick,->-] (5,2) -- (7,2);
	\draw[very thick,->-] (7,2) -- (9,2);
	\draw[very thick,->-] (5,4) -- (7,4);
	\draw[very thick] (7,4) -- (9,4);
	\draw[very thick] (5,0) -- (5,2);
	\draw[very thick,->-] (5,2) -- (5,4);
	\draw[very thick,->-] (7,0) -- (7,2);
	\draw[very thick,->-] (7,2) -- (7,4);
	\draw[very thick,->-] (9,0) -- (9,2);
	\draw[very thick] (9,2) -- (9,4);
	\draw[very thick,->-] (7,0) -- (5,2);
	\draw[very thick,->-] (9,0) -- (7,2);
	\draw[very thick,->-] (7,2) -- (5,4);
	\draw[very thick,->-] (9,2) -- (7,4);
	\node[rectangle,fill=white,gray!80,inner sep=0pt,minimum size=12pt] at (7,2) {}; 
	\node at (7.4,2.4) {$r$};
	
	%% Combined Patch
	\foreach \x in {11,...,13}{
		\foreach \y in {1,...,3}{
			\node[circle,draw,fill=gray!40,inner sep=0pt,minimum size=12pt] at (\x,\y) {}; 
		}
	}
	\node[circle,draw,fill=gray!40,inner sep=0pt,minimum size=12pt] at (12,0) {};
	\node[circle,draw,fill=gray!40,inner sep=0pt,minimum size=12pt] at (12,4) {};
	\node[circle,draw,fill=gray!40,inner sep=0pt,minimum size=12pt] at (10,2) {};
	\node[circle,draw,fill=gray!40,inner sep=0pt,minimum size=12pt] at (14,2) {};
	\node[circle,draw,fill=gray!40,inner sep=0pt,minimum size=12pt] at (14,0) {};
	\node[circle,draw,fill=gray!40,inner sep=0pt,minimum size=12pt] at (13,0) {};
	\node[circle,draw,fill=gray!40,inner sep=0pt,minimum size=12pt] at (14,1) {};
	\node[circle,draw,fill=gray!40,inner sep=0pt,minimum size=12pt] at (10,3) {};
	\node[circle,draw,fill=gray!40,inner sep=0pt,minimum size=12pt] at (10,4) {};
	\node[circle,draw,fill=gray!40,inner sep=0pt,minimum size=12pt] at (11,4) {};
	\draw[very thick] (10,0) -- (12,0);
	\draw[very thick,->-] (12,0) -- (14,0);
	\draw[very thick,->-] (10,2) -- (12,2);
	\draw[very thick,->-] (12,2) -- (14,2);
	\draw[very thick,->-] (10,4) -- (12,4);
	\draw[very thick] (12,4) -- (14,4);
	\draw[very thick] (10,0) -- (10,2);
	\draw[very thick,->-] (10,2) -- (10,4);
	\draw[very thick,->-] (12,0) -- (12,2);
	\draw[very thick,->-] (12,2) -- (12,4);
	\draw[very thick,->-] (14,0) -- (14,2);
	\draw[very thick] (14,2) -- (14,4);
	\draw[very thick,->-] (12,0) -- (10,2);
	\draw[very thick,->-] (14,0) -- (12,2);
	\draw[very thick,->-] (12,2) -- (10,4);
	\draw[very thick,->-] (14,2) -- (12,4);
	\node[rectangle,fill=white,gray!80,inner sep=0pt,minimum size=12pt] at (12,2) {}; 
	\node at (12.4,2.4) {$r$};
	\node[circle,draw,fill=gray!40,inner sep=0pt,minimum size=12pt] at (12,2) {}; 
	\node[circle,fill=black,inner sep=0pt,minimum size=7pt] at (12,2) {}; 
	\node at (11.6,1.6) {$p$};
	
	% Legend
	\node[circle,draw,fill=gray!40,inner sep=0pt,minimum size=12pt] at (0,-1) {};
	\node[right=3mm] at (0,-1) {velocity DoF ($\mathbf{P}_2$)};
	\draw[very thick,->-] (4.5,-1) -- (5.5,-1) {};
	\node[right=1mm] at (5.5,-1) {magnetic field DoF (lowest-order N\'{e}d\'{e}lec)};
	\node[circle,fill=black,inner sep=0pt,minimum size=7pt] at (0,-2) {};
	\node[right=3mm] at (0,-2) {pressure DoF ($P_1$)};
	\node[rectangle,fill=white,gray!80,inner sep=0pt,minimum size=12pt] at (5,-2) {};
	\node[right=3mm] at (5,-2) {Lagrange multiplier DoF ($P_1$)};
	\end{tikzpicture}
	\caption{Groups of DoFs used in the Vanka relaxation schemes.  At left and center are the DoF patches for the fluids and electromagnetic subproblems, respectively, used for the segregated Vanka relaxation scheme.  Both the purist and coupled Vanka schemes use the patches of DoFs depicted at right.  In the purist approach, the same collection of velocity and magnetic-field degrees of freedom are relaxed twice, once with the central nodal pressure DoF included in the Vanka block and once with the electromagnetic Lagrange multiplier DoF.  In the coupled approach, both nodal constraint DoFs at the center of the patch are included in a single Vanka block.}
	\label{fig:Vanka_blocks}
\end{figure}
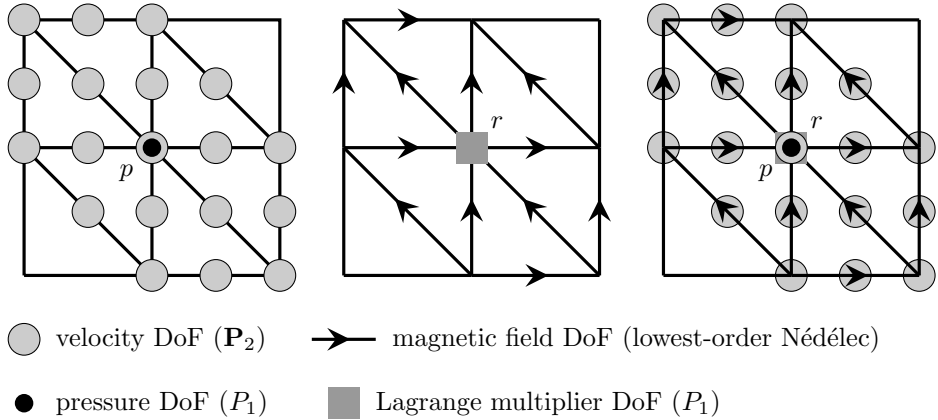

The final variant that we propose is a {\it coupled Vanka} relaxation,
where we take advantage of the fact that the two constraint degrees of
freedom are discretized in the same finite element space and, as such,
have collocated DoFs.  The algebraic viewpoint here is that we loop
over the rows of $\hat{B}$ in pairs, taking one row corresponding to a
pressure DoF and a second row corresponding to a Lagrange multiplier
DoF, both at the same point in the mesh.  This approach naturally
samples the coupling in the $\hat{F}$ matrix, including its
off-diagonal blocks.  From a topological point of view, the coupling is
similarly easy to define, as we construct one Vanka block for each
node in the mesh, composed of the pressure and Lagrange multiplier
DoFs associated with each node along with the velocity and magnetic
field DoFs on the elements adjacent to the node.  These blocks again match those depicted at right of Figure \ref{fig:Vanka_blocks} but, in contrast to purist Vanka, \textit{both} central nodal degrees of freedom are included in a single block, leading to a total number of Vanka blocks equal to the number of vertices in the mesh.

\subsection{Implementation and parallelization}

As a software framework for this work, we use Firedrake
\cite{FRathgeber_etal_2017a} for the finite element discretization
and PETSc \cite{SBalay_etal_2010a, petsc-user-ref} for the
nonlinear and linear solvers.  This pairing was chosen because of the
tight integration between discretizations and solvers in the two
packages \cite{kirby2018solver}, allowing for a natural definition of
the relaxation schemes and multigrid methods detailed above. To drive
the Vanka relaxation, we make use of the PCPatch functionality
introduced to PETSc in \cite{farrell2018augmented,farrell2019c}.
For reproducibility, the codes used to generate the numerical results,
as well as the major Firedrake components needed, have been archived on
Zenodo \cite{zenodo/Firedrake-20200622.1}.

Following the design decision made above, all aspects of the
relaxation scheme are naturally parallelizable.  Because we start with
a coarse triangulation that is refined to create the multigrid
hierarchy, some load-balancing is naturally needed to obtain the best
possible performance. Since the primary workload associated with the
solvers is built around the vertices of the mesh, we choose a parallel
decomposition where we approximately balance the number of
vertices on each processor, independently on every level of the
hierarchy. This is done once when the meshes are initially created
and then used consistently throughout all of the calculations. To account for the need to
compute residuals for each DoF within each patch, a nodal-distance-2
halo is needed for all operations within the relaxation scheme to be
computable without additional communications, and the parallel mesh
distribution is done consistently with this \cite{MLange_etal_2016a}.

\section{Numerical Results}
\label{sec:numerical_results}

In this section, we present numerical results for three test
problems: steady-state Hartmann flow in two dimensions; a time-dependent
two-dimensional island coalescence simulation; and a three-dimensional
(steady) MHD generator test problem.  All timings reported were
measured on the Compute Canada machine \texttt{Niagara}, a cluster of
1548 Lenovo SD530 servers, each with 40 2.4 GHz Intel Skylake cores,
and 202 GB of RAM.

\subsection{Two-Dimensional Hartmann Flow}
\label{ssec:hartmann}

First, we consider a steady-state Hartmann flow model posed over the square
domain $\Omega = [-1/2,1/2]^2$.  This model problem considers a
section of a duct or channel through which a conducting fluid is
flowing and is subjected to a transverse magnetic field, $\vec{B}_0 =
(0,B_0)^{\top}$, applied in the perpendicular direction to the flow.  The
flow itself is driven in the positive $x$-direction by an applied
pressure gradient, $\frac{\partial p}{\partial x} = -G_0$.  The
analytical solution to this system is given by vector fields $\vec{u}
= (u_1(y),0)^{\top}$ and $\vec{B} = (B_1(y),B_0)^{\top}$, with
\begin{align*}
u_1(y) & = \frac{GRe}{2\Ha\tanh(\Ha/2)}\left(1-\frac{\cosh(y\Ha)}{\cosh(\Ha/2)}\right), \\
B_1(y) & = \frac{G}{2}\left(\frac{\sinh(y\Ha)}{\sinh(\Ha/2)} -
2y\right),
\end{align*}
where $Re$ and $Re_m$ are prescribed fluid and magnetic Reynolds
numbers, $\Ha = \sqrt{ReRe_m}$ is the so-called Hartmann number, and $G = \frac{2\Ha
  \sinh(\Ha/2)}{Re(\cosh(\Ha/2)-1)}$.  The fluid pressure is given by
$p(x,y) = -G x - B_1^2(y)/2$ (plus any scalar value), and the Lagrange
multiplier $r=0$.  We apply Dirichlet boundary conditions to both
components of $\vec{u}$ and $\vec{B}$, as well as to $r$, on all
edges.  To eliminate the trivial pressure nullspace, we fix the value
of the pressure at the origin to be 0.

We consider a weak-scaling experiment, fixing the finest grid to be a
$120\times 120$ mesh for serial runs (189K total degrees of freedom),
$480\times 480$ for runs on 16 cores (3.00M degrees of freedom), and
$1920\times 1920$ for runs on 256 cores (47.95M degrees of
freedom). All experiments use a fixed $15\times 15$ coarsest grid,
resulting in 4 levels in the multigrid hierarchy for the serial
runs, 6 levels for the runs on 16 cores, and 8 levels for the runs on
256 cores.  For these runs, we used a single node of Niagara for
single- and 16-core runs, and 32 cores on each of 8 nodes for 256-core
runs.  Overall, we use a multigrid V(2,2) cycle as a preconditioner
for FGMRES, stopping when either the relative reduction in the norm of
the residual reaches $10^{-6}$ or when the absolute value of the norm
of the residual reaches $10^{-6}$ (whichever occurs first).  The pairs
of pre- and post-relaxation sweeps are combined via Chebyshev
acceleration, with hand-tuned parameters for the interval over which
the relaxation is most effective.  That is, $k$ sweeps of Vanka relaxation correspond to $k$ steps of the form,
\[ u_{i+1} = u_i + \omega_i M^{-1} \left ( b - A u_i\right ),\]
where $M^{-1}$ corresponds to the combined action of all $N$ additive Vanka block updates and the $\omega_i$'s are chosen to minimize a Chebyshev polynomial over an interval $[i_1,i_2]$.  The interval choice, which depends on an estimate of the spectral radius of $M^{-1}A$ and the multigrid coarsening rate, can be estimated using PETSc.  However, we use hand-tuned intervals to avoid some irregularities of the estimated values. These intervals are $[1.5,8.0]$ for the segregated method, 
$[1.5,16.0]$ for purist Vanka, and $[2.0,8.0]$ for the coupled approach.
All other parameters of the Vanka methods are kept
constant between runs, aside from the patch constructions detailed
above.  The coarsest grid is solved using a direct solve (via MUMPS \cite{amestoy2001}).

Table \ref{tab:Hartmann_itns} shows the averaged number of
multigrid-preconditioned FGMRES iterations per linear solve along with
the number of Newton steps required to reach solution, prescribed as a
relative reduction in the nonlinear residual norm by a factor of
$10^5$.  As expected, we see growth in the effort needed to solve
these problems with increasing values of $Re$ and $Re_m$, reflected in
both the number of Newton iterations required and number of
preconditioned FGMRES iterations needed to reach convergence for the
linear solves.  Note, however, that there is some decrease in the
number of required Newton steps with increasing problem size at fixed
parameters.  This is also to be expected as the relative scaling
between the advective and diffusive terms in the discretized systems
becomes somewhat more favorable as the grid is refined.

\begin{table}
  \centering
  \caption{Averaged multigrid-preconditioned FGMRES iteration counts and number of nonlinear
    iterations (subscripted) for Hartmann test problem in 2D. Methods are Segregated (seg), Purist (pur), and Coupled (coup).}
  \label{tab:Hartmann_itns}
  \begin{tabular}{c||c@{\hskip 5pt}c@{\hskip 5pt}c|c@{\hskip 5pt}c@{\hskip 5pt}c|c@{\hskip 5pt}c@{\hskip 5pt}c}
    \toprule
    $\!\!\!(Re,Re_m)\!\!\!$& \multicolumn{3}{c|}{$120\times 120$}
    & \multicolumn{3}{c|}{$480\times 480$}
    & \multicolumn{3}{c}{$1920\times 1920$}\\
    & seg & pur & coup     & seg & pur & coup     & seg & pur & coup \\
    \midrule
    $(4,4)$ &   $6.33_3$ & $8.33_3$ & $6.67_3$ & $5.00_2$ & $10.50_2$ & $5.50_2$ & $5.00_2$ & $10.00_2$ & $5.50_2$ \\
    $(16,4)$ &  $6.67_3$ & $9.67_3$ & $6.67_3$ & $7.50_2$ & $10.50_2$ & $8.00_2$ & $5.00_2$ & $10.00_2$ & $5.50_2$ \\
    $(64,4)$ &  $8.67_3$ & $11.00_3$ & $8.33_3$ & $10.00_3$ & $13.00_3$ & $10.33_3$ & $5.50_2$ & $9.50_2$ & $6.00_2$ \\
    $(4,16)$ &  $9.33_3$ & $11.33_3$ & $9.67_3$ & $9.33_3$ & $11.00_3$ & $10.33_3$ & $8.67_3$ & $14.67_3$ & $9.67_3$ \\
    $(16,16)$ & $8.25_4$ & $10.00_4$ & $8.50_4$ & $9.67_3$ & $13.33_3$ & $10.67_3$ & $9.00_3$ & $15.00_3$ & $10.00_3$ \\
    $(64,16)$ & $8.75_4$ & $11.00_4$ & $9.25_4$ & $11.33_3$ & $14.00_3$ & $12.00_3$ & $10.33_3$ & $17.67_3$ & $11.00_3$ \\
    $(4,64)$ &  $11.60_5$ & $14.40_5$ & $11.60_5$ & $12.50_4$ & $15.75_4$ & $11.25_4$ & $14.00_4$ & $19.50_4$ & $14.00_4$ \\
    $(16,64)$ & $12.60_5$ & $15.20_5$ & $13.20_5$ & $14.00_5$ & $15.40_5$ & $12.20_5$ & $14.00_4$ & $23.50_4$ & $12.50_4$ \\
    $(64,64)$ & $28.17_6$ & $16.33_6$ & $13.50_6$ & $17.60_5$ & $17.80_5$ & $15.00_5$ & $25.75_4$ & $30.00_4$ & $15.00_4$ \\
    \bottomrule
    \end{tabular}
  \end{table}

Comparing approaches to the relaxation, Table \ref{tab:Hartmann_itns}
clearly shows the benefit of the coupled Vanka scheme.  For the segregated
method, we see clear increases in iteration counts per linear solve
as both $Re$ and $Re_m$ increase.  This is somewhat expected, since
this approach neglects the coupling between $\vec{u}$ and $\vec{B}$ in
the relaxation scheme, but this coupling becomes more important as
$\Ha$ increases.  On the finest grid, we see a factor of 5 increase
in averaged iteration counts per linear solve between the lowest and
largest $\Ha$ values.  This variation with $\Ha$ is ameliorated
in the results for the purist and coupled approaches.  For the purist
approach, we see an increase of a factor of 2-3 between the minimum
and maximum values of the average number of linear iterations per
Newton step for each mesh size.  While iteration counts are quite
steady for each set of parameter values between the $120\times 120$
and $480\times 480$ mesh, we see some growth in iteration counts on
the $1920\times 1920$ mesh, particularly for larger values of $Re_m$.
In contrast, the coupled approach is similarly scalable with parameter
values, showing an increase in average number of linear iterations per
Newton step by a factor of less than 3 for all mesh sizes, but with
near-perfect scalability across mesh sizes at fixed parameter values.

Much of the variation in linear and nonlinear solver performance is
directly reflected in wall-clock time-to-solution, as shown in Table
\ref{tab:Hartmann_time}.  In particular, for each method on any single
mesh, we see the CPU-time scaling is dominated by the total number of
FGMRES iterations required for convergence.  Comparing between
methods, we see the high relative cost per iteration of the purist
method that, as expected, costs about twice as much per iteration as
the segregated or coupled approaches.  

\begin{table}
  \centering
  \caption{Solution time (in minutes) for Hartmann test problem in 2D.  CPU
    times are in serial for $120\times 120$ mesh, on 16 cores for
    $480\times 480$ mesh, on 256 cores for $1920\times 1920$ mesh.
    Methods are Segregated (seg), Purist (pur), and Coupled (coup).}
  \label{tab:Hartmann_time}
  \begin{tabular}{c||ccc|ccc|ccc}
    \toprule
    $(Re,Re_m)$& \multicolumn{3}{c|}{$120\times 120$}
    & \multicolumn{3}{c|}{$480\times 480$}
    & \multicolumn{3}{c}{$1920\times 1920$}\\
    & seg & pur & coup     & seg & pur & coup     & seg & pur & coup \\
    \midrule
    $(4,4)$ &   0.52 & 0.85 & 0.42 & 0.57 & 1.06 & 0.43 & 0.75 & 1.32 & 0.55 \\
    $(16,4)$ &  0.54 & 0.92 & 0.42 & 0.64 & 1.01 & 0.49 & 0.75 & 1.32 & 0.55 \\
    $(64,4)$ &  0.60 & 0.99 & 0.47 & 0.89 & 1.50 & 0.74 & 0.77 & 1.30 & 0.57 \\
    $(4,16)$ &  0.62 & 1.00 & 0.51 & 0.87 & 1.38 & 0.71 & 1.09 & 2.15 & 0.93 \\
    $(16,16)$ & 0.69 & 1.14 & 0.58 & 0.88 & 1.55 & 0.75 & 1.11 & 2.15 & 0.94 \\
    $(64,16)$ & 0.71 & 1.21 & 0.61 & 0.95 & 1.59 & 0.78 & 1.18 & 2.44 & 0.98 \\
    $(4,64)$ &  0.98 & 1.72 & 0.84 & 1.16 & 2.06 & 0.96 & 1.69 & 3.33 & 1.44 \\
    $(16,64)$ & 1.03 & 1.79 & 0.93 & 1.53 & 2.49 & 1.19 & 1.73 & 3.86 & 1.32 \\
    $(64,64)$ & 2.22 & 2.20 & 1.10 & 1.76 & 2.77 & 1.41 & 2.66 & 4.58 & 1.53 \\
    \bottomrule
    \end{tabular}
  \end{table}

Measuring total memory usage (as reported by Slurm's \texttt{sacct}),
we see that the serial runs of both coupled and segregated
relaxation schemes peak at about 2GB of usage while purist relaxation
requires about 3GB of memory, with only slight variation with
parameters.  In comparison, a direct solve (using MUMPS) of this
problem required about 1GB of memory, so this represents some notable
overhead for a small 2D mesh, but we note that all aspects of the
multigrid solvers (except the coarsest grid solve) are naturally
optimally parallelizable.  For the $480\times 480$ mesh on 16 cores,
peak memory usage over all cores was similar, rising to about 3.5GB
for purist relaxation, while staying at about 2GB for both coupled and
segregated.  Slightly more growth was seen for the $1920\times 1920$
mesh on 256 cores, with peak memory usage over all cores of about
2.5GB for both coupled and segregated relaxation schemes, and 4GB per
core for purist relaxation.  This may be due to effects of the
parallel partitioning not being perfectly balanced between threads;
however, this still shows the expected near-optimal scaling of memory
consumption of multigrid algorithms.  We note that the extra overhead in
memory usage by the purist relaxation is to be expected, since it
creates twice as many patch matrices as the coupled algorithm, but
each is roughly the same size as those used in the coupled approach.
Similarly, the purist algorithm creates an equal number of patches to
the segregated approach, but each is about twice as large as the
patches in the segregated approach.  Thus, the memory consumption due
solely to the relaxation schemes is about double for the purist
approach than for either the segregated or coupled approaches.

Overall, the results in Tables \ref{tab:Hartmann_itns} and
\ref{tab:Hartmann_time} highlight the effectiveness of the coupled
relaxation strategy, consistently yielding the lowest CPU times to
solution of all the approaches tested herein and the lowest memory
usage.  As discussed above, this is not terribly surprising.  The
segregated approach has a similar memory overhead as the coupled
approach, but neglects important coupling between the fluid and
magnetics degrees of freedom, so is not expected to be as robust as
the coupled approach.  In contrast, the purist approach retains that
coupling, but we see little benefit to its much higher memory
footprint and cost-per-iteration.  Consequently, in the rest of this
paper, we present results only for the coupled relaxation scheme.

To effectively solve problems at larger Hartmann numbers, two
ingredients are needed.  First of all, as the physical parameters
increase, the convergence of Newton's method degrades, as seen in
Table \ref{tab:Hartmann_itns}; for large-enough parameters
(not shown here), Newton's method fails to converge at all.  To
address this, we use a continuation approach, where we use the
solution at low Hartmann number as the initial guess for Newton's
method at larger parameter values.  As shown below, we can use such a
strategy to solve problems at much larger parameter values than
possible using a fixed initial guess.  The
second necessary ingredient is the use of \textit{shallower} multigrid
hierarchies, with finer grids chosen as the coarsest mesh.  To a large
extent, this is intuitive: the discretization that we consider has no
stabilization to deal with large convective terms in either the
momentum or electromagnetic equations and, thus, is expected to break
down as those terms become dominant.  If the quality of discretization
on the coarsest grid becomes sufficiently poor, then it will provide a
poor correction to the finest grid discretization, even if that
discretization is relatively faithful to the continuum problem.  Thus,
as we seek solutions at increasingly large Hartmann numbers, we must
use increasingly large coarsest grids.

\begin{figure}
  \begin{center}
    \includegraphics{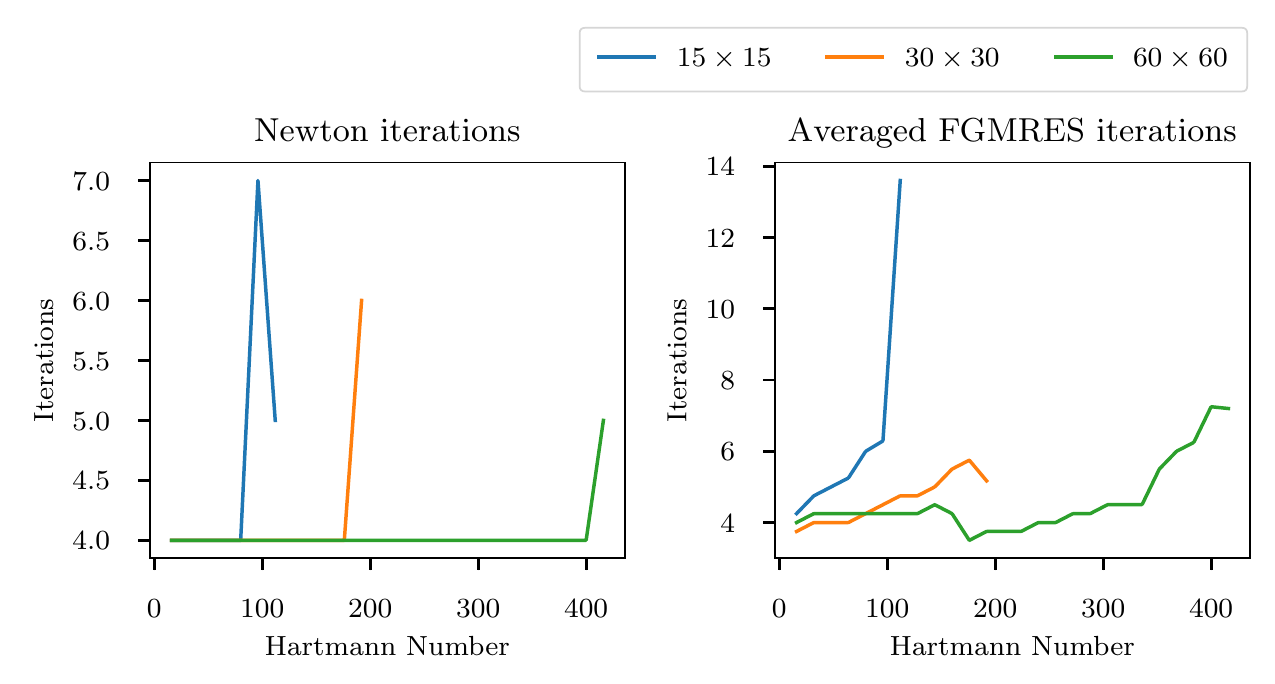}
    \caption{Number of Newton steps (left) and average number of
      iterations of preconditioned FGMRES using coupled Vanka to solve each linearization
      (right) when using continuation to solve the Hartmann problem
      with $Re = Re_m = \Ha$ on the $480\times 480$ mesh, with
      coarsest grids in the hierarchy of $15\times 15$, $30\times 30$,
      and $60\times 60$.}
    \label{fig:Hartmann_continuation}
	\end{center}
\end{figure}

To illustrate this behavior,
Figure \ref{fig:Hartmann_continuation} shows nonlinear and averaged
linear iteration counts for solution of the problem with $Re = Re_m =
\Ha$ on the $480\times 480$ mesh, using continuation in steps of $16$
from $Re = Re_m = \Ha = 16$ with $15\times 15$, $30\times 30$, and
$60\times 60$ meshes as the coarsest grid (leading to 6, 5, and 4
levels in the multigrid hierarchy, respectively).  In these results,
we use the Eisenstat-Walker inexact Newton approach
\cite{SCEisenstat_HFWalker_1996a} to set the stopping criteria for the
linear solver, with the same nonlinear solver tolerances as above.  As
we see, performance of Newton's method is quite stable until
convergence fails, and the same is true for that of the linear solver.  Moreover, when failure occurs, it is the linear solver that breaks down first, giving stalling convergence on the first Newton step at a new value of $Re = Re_m = \Ha$.
What is, perhaps, most impressive in these results is just how
large the problem parameters can become before convergence fails: with
a $15\times 15$ coarsest mesh, continuation successfully finds
solutions past $\Ha = 100$, and this parameter value doubles with each
refinement of the coarsest mesh.

\subsection{Two-dimensional island coalescence}
\label{ssec:island}

\begin{figure}
	\begin{center}
		\includegraphics[width=0.6\linewidth]{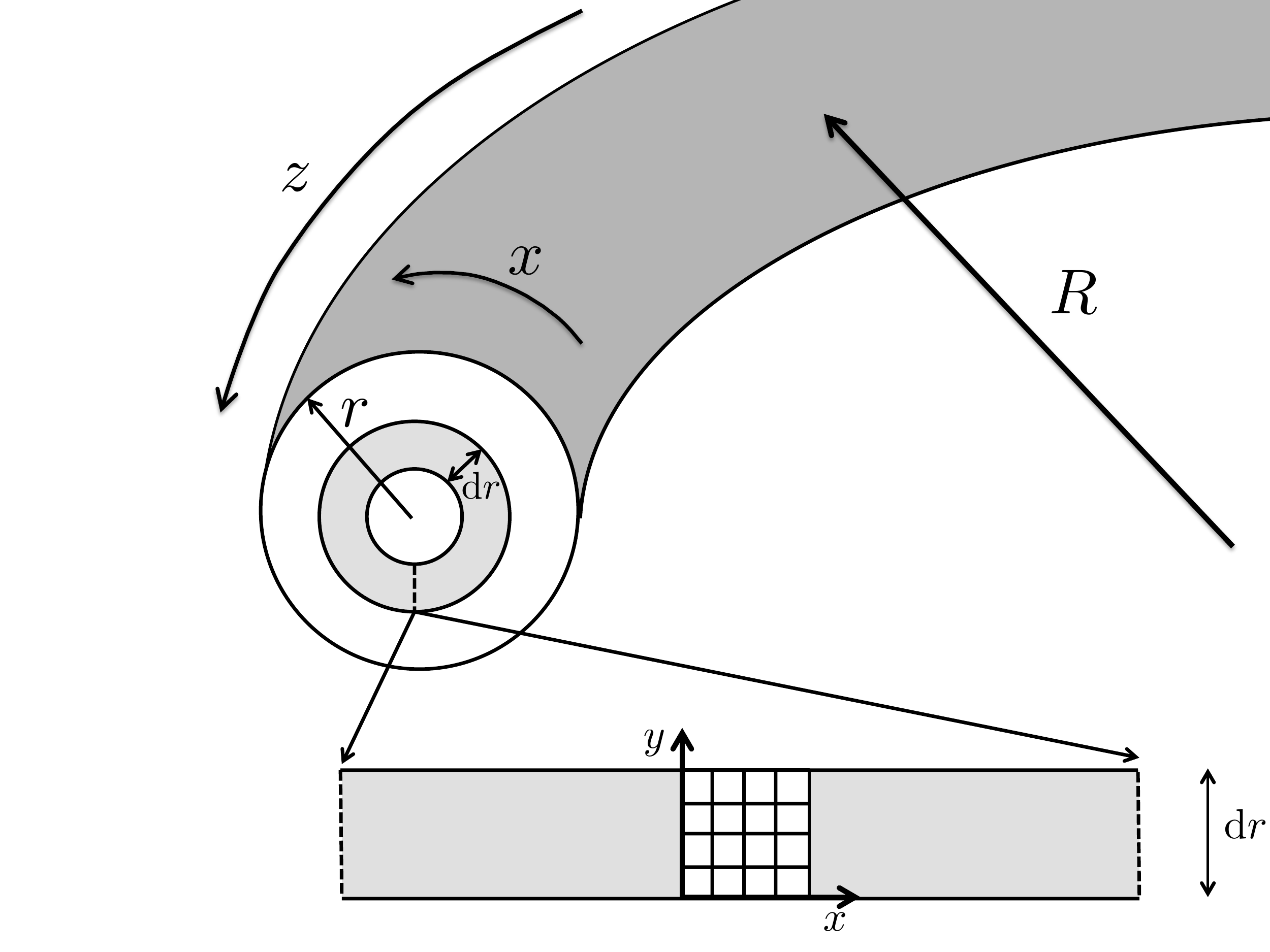}
		\caption{Cross-sectional view of large aspect-ratio tokamak
			geometry, with major radius, $R$, and minor radius, $r$,
			satisfying $R \gg r$.  A cross-section of thickness d$r$ can be
			unfolded to create a Cartesian grid as pictured.}
		\label{fig:islanddomain}
	\end{center}
\end{figure}

We next consider a time-\\dependent two-dimensional test problem that
mimics magnetic reconnection in a large aspect ratio tokamak
\cite{1978BatemanG-aa,1976StraussH-aa}.  In this setting, we consider the
cross-section of flow of a magnetically confined plasma, in which a
large magnetic field is imposed in the toroidal direction, resulting
in effectively two-dimensional dynamics.  We then consider an annulus
in the cross-sectional direction that is unfolded and rescaled to
make a square domain, $\Omega = [-1,1]^2$, with a periodic mapping
between its right and left edges, see Figure \ref{fig:islanddomain}.

On this domain, we consider a model
problem of a physical instability that can arise from perturbations of
the magnetic current density, known as island coalescence, which has been simulated in numerous works (eg. \cite{ECCyr_EtAl_2013a,2008ChaconL-aa,LChacon_DAKnoll_2006a,2002ChaconL_KnollD_FinnJ-aa,2008PhilipB_ChaconL_PerniceM-aa,JAdler_EtAl_2013a}).  In this
problem, two initially isolated ``islands'' of current density are
perturbed, resulting in a breaking and then \emph{reconnection} of the magnetic field lines. At the reconnection point, where the magnetic field comes back together, a sharp peak in current density
occurs. 

\begin{figure}
  \begin{center}
    \includegraphics{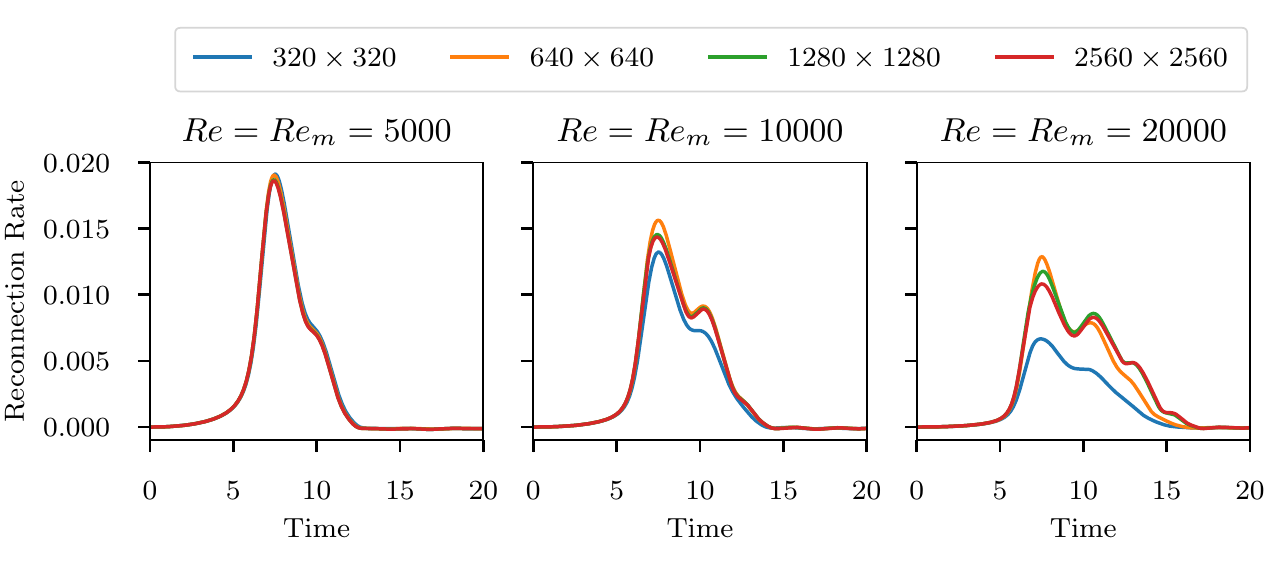}
    \caption{Reconnection rates computed for island coalescence model
      with $Re = Re_m$ values of 5,000 (left), 10,000 (middle),
      and 20,000 (right).}
    \label{fig:reconnection}
  \end{center}
\end{figure}

As an equilibrium solution, we take
\begin{align*}
  \vec{u}_0(x,y) & = \vec{0}, \\
  \vec{B}_0(x,y) & = \frac{1}{\cosh(2\pi y) +
k\cos(2\pi x)}\left(\begin{array}{c} \sinh(2\pi y) \\ k\sin(2\pi
    x)\end{array}\right), \\
  p(x,y) & = \frac{1-k^2}{2}\left(1+\frac{1}{\left(\cosh(2\pi y) +
    k\cos(2\pi x)\right)^2}\right), \\
  r(x,y) & = 0,
\end{align*}
for $k=0.2$.  To support this as an equilibrium solution of Equations
\eqref{eqn:ns}--\eqref{eqn:div-b-free} requires imposition of
right-hand side functions $\vec{f} = \vec{0}$ and
\[
\vec{g} = \frac{-8\pi^2(k^2-1)}{Re_m\left(\cosh(2\pi y) + k\cos(2\pi
  x)\right)^3}\left(\begin{array}{c} \sinh(2\pi y) \\ k\sin(2\pi
  x)\end{array}\right).
\]
As an initial condition, we perturb the equilibrium solution for the
magnetic field by
\[
\delta\vec{B} = \frac{\epsilon}{\pi}
\left(\begin{array}{c} -\cos(\pi x)\sin(\pi y/2) \\
  \cos(\pi y/2)\sin(\pi x)/2\end{array}\right),
\]
taking $\epsilon = -0.01$.  We measure physical fidelity of the
simulation by considering the \textit{reconnection rate} of the
simulation, which is the time rate of change of the poloidal flux function, $\Psi$, such that $\vec{B} = -\nabla\times\left (0,0,\Psi\right )^{\top}$.  In our formulation, this can be computed as the difference between $\nabla\times \vec{B}$
at the reconnection point (origin for our computational domain) at the current time to its initial value, scaled by
$1/\sqrt{Re_m}$.  For more details, see \cite{1978BatemanG-aa,1976StraussH-aa}.  For low magnetic Reynolds numbers, the area of reconnection is wider with a less steep gradient in the current density when
the peak occurs. As the magnetic Reynolds number increases, this reconnection zone narrows,
resulting in a sharper, yet shorter, peak. In addition, a ``sloshing" effect occurs, where
the islands bounce a little before fully merging into one \cite{LChacon_DAKnoll_2006a}. This yields a peak in the
reconnection rate, whose height oscillates as the islands come together. More sloshing occurs for higher magnetic Reynolds numbers.  This is validated in our computational tests shown in Figure \ref{fig:reconnection}.

While the curl of a N\'ed\'elec field is naturally
computed in a discontinuous Galerkin space (via the usual exact
sequence), because we are interested in a pointwise measurement of the
quantity, we instead project $\nabla\times\vec{B}$ onto the continuous
space of piecewise linear functions for this measurement.  This requires a solve with the mass matrix for the piecewise linear space, which we do using a simple multigrid method as a preconditioner for the conjugate gradient iteration.  The preconditioner is formed once before time-stepping begins, and reused in a postprocessing step after each time step.  We do not time this postprocessing step in the results below, so we have not made significant effort to optimize this computation.   Since
symmetry around the origin is important in this measurement, we
discretize this problem on ``crossed'' triangular meshes, taking
uniform quadrilateral meshes of $[-1,1]^2$ and cutting each square
cell into 4 triangles, adding a vertex at the cell midpoint.
  
We discretize the temporal derivatives using the L-stable BDF2 method.
For each spatial mesh, a fixed time-step is used; however, to start
the time-stepping, a refinement of the first time-step into ten
(equal) substeps is used, as we observed nonlinear convergence issues
without this.  For the first substep, a Crank-Nicolson discretization
is used (to avoid the need for additional initial data for BDF2).
There are two relevant CFL numbers for these simulations.  The
classical fluid CFL, which we take to be $u_{\text{max}} \frac{\Delta t}{h}$, where
$u_{\text{max}}$ is the maximum magnitude of the velocity vector,
$\vec{u}$, at a given time step, $\Delta t$ is the time step itself, and
$h$ is the spatial mesh width, calculated as the length of the
shortest edges of the crossed triangular mesh.  Similarly, we consider
the Alfv\'{e}n CFL, defined as $B_{\text{max}} \frac{\Delta t}{h}$, where $B_{\text
  max}$ is the maximum magnitude of the magnetic field, $\vec{B}$.
For these calculations, we approximate $u_{\text{max}}$ and
$B_{\text{max}}$ by projecting $\vec{u}\cdot\vec{u}$ and
$\vec{B}\cdot\vec{B}$ onto the scalar piecewise-constant space and
taking the square root of the maximum elementwise value of that projection.

We consider this problem with $Re = Re_m$ values of 5,000, 10,000,
and 20,000, on four levels of spatial refinement, from $320\times 320$
quadrilateral elements, each cut into 4 triangles in the ``crossed''
mesh, up to $2560\times 2560$ cut quadrilateral elements.  For the
finest mesh, the discretization leads to a discrete system with about
170 million degrees of freedom.  For all meshes, we consider a
multigrid hierarchy formed from a coarsest grid of size $20\times 20$.  As we use a
second-order time-stepping scheme, we consider a time-step of $\Delta t =
0.025$ on the coarsest level, and halve the time-step with each
(uniform) refinement of the spatial mesh.  Figure
\ref{fig:reconnection} shows the reconnection rates for these
problems, computed across the meshes.  For the $\Re_m = 5,000$ case, we
see that the reconnection rate is close to convergence already on the
coarsest mesh, but that finer meshes are needed to accurately resolve
the true dynamics at higher $\Re_m$ values.  For runs with $\Re_m = 40000$
(not shown here), convergence of the reconnection rate was not yet
observed by the $2560\times 2560$ mesh.  Additionally, we see the expected qualitative behavior: a reduction in reconnection rate and increase in ``sloshing" with $Re_m$.  

\begin{figure}
  \begin{center}
    \includegraphics{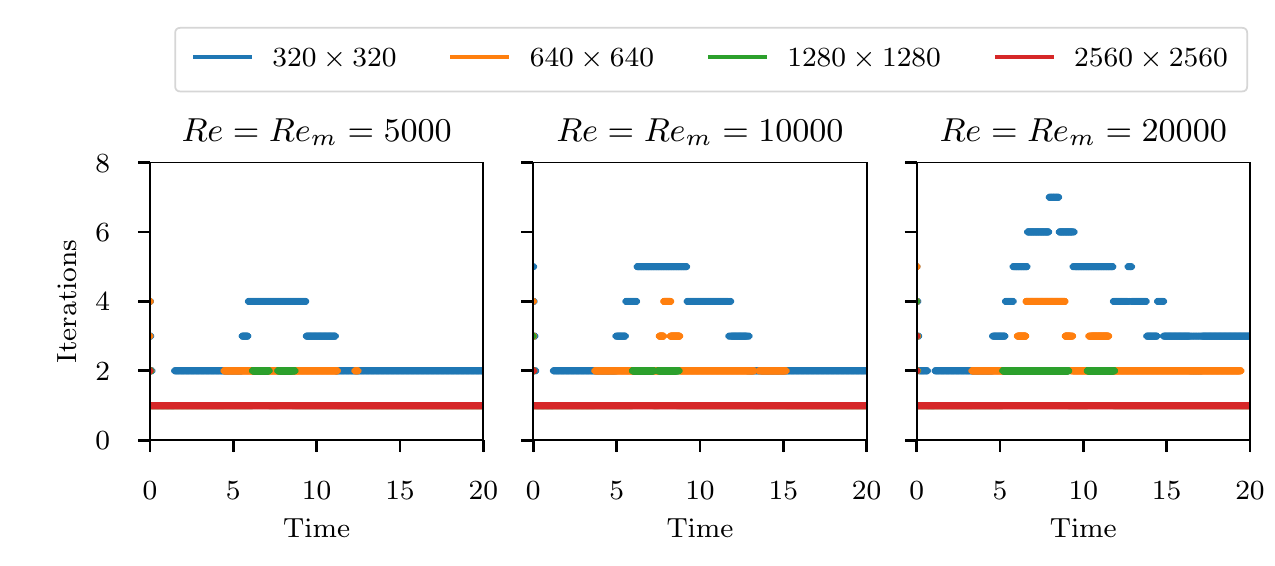}
    \caption{Number of linear solver iterations needed per time-step
      for island coalescence model with $Re = Re_m$ values of
      5,000 (left), 10,000 (middle), and 20,000 (right).}
    \label{fig:island_iterations}
  \end{center}
\end{figure}

\begin{figure}
  \begin{center}
    \includegraphics{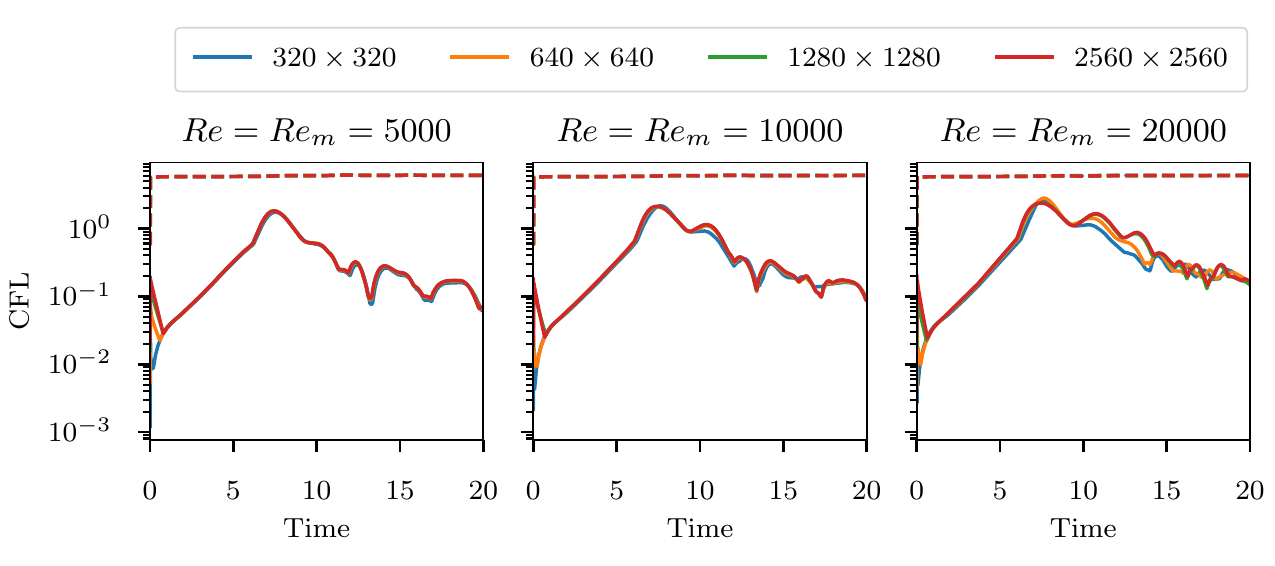}
    \caption{Measured CFL values at each time-step
      for island coalescence model with $Re = Re_m$ values of
      5,000 (left), 10,000 (middle), and 20,000 (right).  Solid lines
      denote the fluid CFL, dashed lines denote the Alfv\'{e}n CFL.}
    \label{fig:island_CFL}
  \end{center}
\end{figure}

Figure \ref{fig:island_iterations} shows the number of linear
iterations needed per time-step for these simulations, again only
considering the coupled Vanka relaxation scheme.  Here, based on
preliminary experiments, we used V(3,3) multigrid cycles with the Chebyshev
interval chosen to be $[2.0,10.0]$.  For each
time-step, we use the solution from the previous time-step as the
initial guess for Newton's method applied to the discretized nonlinear
equations at the new time-step.  Convergence is measured by requiring
the $\ell_2$ norm of the nonlinear residual be reduced by either a
relative factor of $10^{8}$ or to an absolute value below $10^{-6}$.
Since time-stepping provides a good-quality initial guess, we observe
that we typically require only one or two Newton steps to achieve
convergence based on the absolute residual-norm tolerance.  For
coarser meshes, we observe some growth in the required number of
linear iterations per time step as we increase $Re_m$, but this is
largely mitigated on the finest meshes, where only a single linear
iteration is needed per time step throughout the simulation.  We note
that the time steps used are chosen so that the Alfv\'{e}n CFL number is
larger than one (except in the first time step, where sub-steps are
used to initialize the simulation), and kept at a constant factor of
about 6 over all simulations, as shown in Figure
\ref{fig:island_CFL}.  The fluid CFL is much more variable, peaking at
values around 2 at the point where peak reconnection is observed.  We
note that the discretization has no stabilization to cope with the
advective terms, and that problems with nonlinear solver convergence
arose when simulations (not shown here) were run with higher fluid CFL
values, as might be expected.

\begin{figure}
  \begin{center}
    \includegraphics{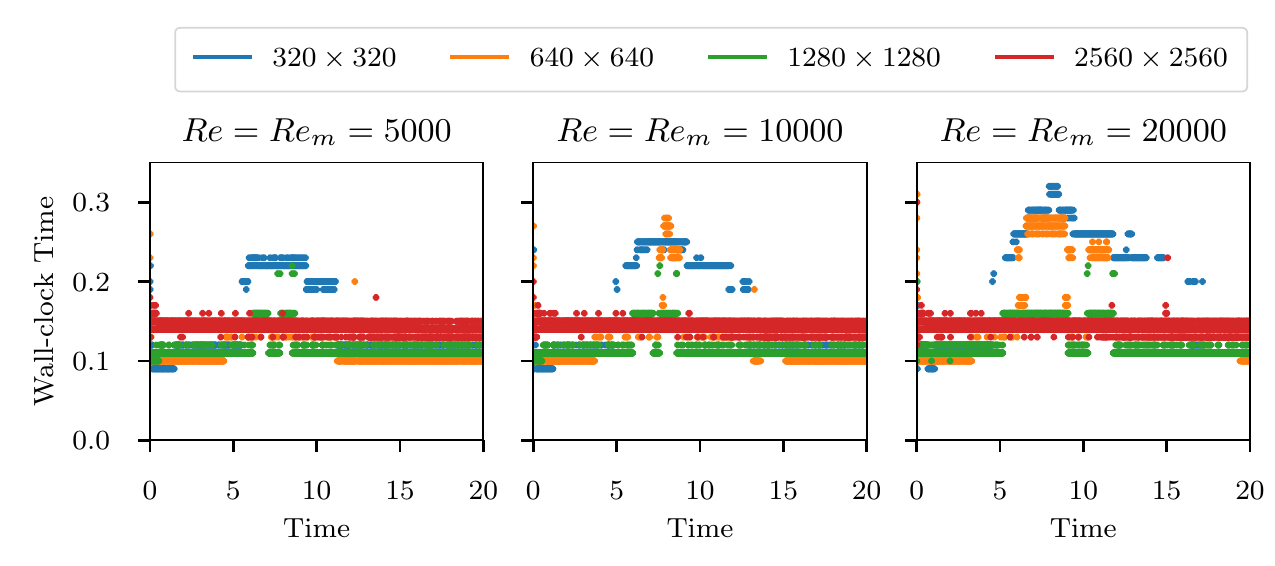}
    \caption{Wall-clock solution time needed per time-step
      for island coalescence model with $Re = Re_m$ values of
      5,000 (left), 10,000 (middle), and 20,000 (right).}
    \label{fig:island_runtimes}
  \end{center}
\end{figure}

A standard weak scaling experiment was performed.
The simulations on the $320\times 320$ mesh were
run on 10 cores of a single node on Niagara, requiring roughly 30GB of
total memory and between 2 and 2.5 hours of wall-clock runtime,
including initialization, input/output, and some postprocessing (such
as computing the CFL values at each time-step, as shown above).  With
each refinement, the number of cores allotted was increased by a
factor of four, leading to 40 cores (1 node) for the $640\times 640$
mesh, 160 cores (4 nodes) for the $1280\times 1280$ mesh, and 640
nodes (16 cores) for the $2560\times 2560$ mesh.  Reported memory
usage for all of these simulations peaks at about 110GB of memory per
node.  Figure \ref{fig:island_runtimes} shows wall-clock runtimes for
the nonlinear system solve at each timestep (including the linear system solve(s) required).
Comparing with Figure \ref{fig:island_iterations}, we see that the
time per time step scales roughly with the number of linear iterations
(as expected), and that there is not undue growth with problem size.
Thus, the total wall-clock runtimes for these simulations grows
roughly by a factor of 2 per refinement, due to the doubling of the
number of time steps with each mesh refinement.

With such low iteration counts, a natural question to ask is whether or not the problem could be efficiently solved using a simpler preconditioner.  To test this hypothesis, we considered the same simulations using just the multigrid relaxation scheme as the preconditioner, with no coarse-grid correction.  On the $640\times 640$ mesh, at $Re_m = 20000$, the initial timesteps require substantially more work, requiring up to 200 linear iterations per timestep, taking over 1 minute per timestep to compute, again on 40 cores.  After some ``spinup'' time, the simulation becomes less expensive, still requiring about 40 linear iterations per timestep, costing about 0.3 seconds per timestep, around three times the cost per timestep of the multigrid simulation.  Similar results were seen at lower Reynolds numbers, although somewhat mitigated, for example requiring 25-30 linear iterations per timestep at $Re_m = 5000$, but still with a cost above double that of the multigrid simulation.  While these results improve somewhat with mesh refinement, the need for and value of the coarse-grid correction process are clear.

\subsection{MHD Generator}

Finally, we consider a three-dimensional steady-state test problem,
where a duct flow is used to generate an electric current.  Similar
simulations have been considered in \cite{shadid2016scalable,
  EGPhillips_etal_2016a, MWathen_CGreif_2020a}.  We use similar
problem specifications as in \cite{MWathen_CGreif_2020a}, considering
the domain $[0,5]\times[0,1]\times[0,1]$.  For the velocity, we
impose an inflow Dirichlet boundary condition of
$\vec{u} = (1,0,0)^{\top}$ on the face with $x=0$, a natural boundary
condition on the face with $x=5$, and no-slip boundary conditions of
$\vec{u}=(0,0,0)^{\top}$ on the four other faces of the domain.  On
the ``bottom'' and ``top'' faces of the domain (with $z=0$ and $z=1$),
we impose a Dirichlet boundary condition on the tangential component
of $\vec{B}$, that is $\vec{n}\times\vec{B} = \vec{n}\times
(0,b_y(x),0)^{\top}$ for outward unit normal vector $\vec{n}$, and
\[
b_y(x) = \frac{B_0}{2}\left(\tanh\left(\frac{x-x_{\rm
    on}}{\delta}\right) - \tanh\left(\frac{x-x_{\rm
    off}}{\delta}\right)\right),
\]
where $B_0 = 1$, $\delta = 0.1$, $x_{\rm on} = 2.0$, and $x_{\rm off}
= 2.5$.  We fix the values of both the pressure and the Lagrange
multiplier to be zero at the origin but, otherwise, impose no boundary
conditions on these.  For these simulations, we fix $Re = Re_m = \Ha =
10$, as in the experiments of \cite{MWathen_CGreif_2020a}.

Table \ref{tab:generator} presents linear and nonlinear convergence
results for the test problem at 3 different levels of resolution using
the coupled Vanka relaxation scheme.  For this example, we use V(2,2)
multigrid cycles with the Chebyshev relaxation interval chosen to be $[2.25,18.0]$.
These meshes are created by constructing the coarsest mesh in the
multigrid hierarchy, first as a mesh of hexahedra then ``cutting''
each hexahedron into six tetrahedra in the usual manner.  These
coarsest meshes are then refined uniformly as needed to generate the
finest meshes in the hierarchy.  Thus, the 2-level solver for the
smallest problem size, the 3-level solver for the middle problem size,
and the 4-level solver for the largest problem size all use the same
coarsest mesh, starting from a $20\times 4\times 4$ mesh of hexahedra
on the domain $[0,5]\times[0,1]\times[0,1]$ and refining one to three
times.  A consequence of this is that the different multigrid solvers
for the same finest-grid problem sizes are solving slightly different
problems, as the processes of cutting hexahedra into tetrahedra and
regular refinement do not perfectly commute for these meshes.
Nonetheless, we believe these problems present a fair comparison of
the solver performance as we add levels to the multigrid hierarchies.
For the nonlinear solver, we again use Newton's method with
Eisenstat-Walker for the stopping criteria.  As expected, we see that the
performance of Newton's method is largely insensitive to the details
of the linear solver.  For both of the largest two problem sizes, we
see some notable degradation in multigrid solver performance with the
4-level hierarchies, which we again ascribe to the poor resolution of
this problem on the coarsest meshes of the hierarchy.  For the middle
problem size with the 4-level hierarchy, the underlying hexahedral
cells are cubes with sidelength of $0.5$, while they are cubes of
sidelength $0.25$ for the 4-level hierarchy and the largest problem
size.  In both cases, these significantly under-resolve the Reynolds
numbers, so it is unsurprising that solver performance is compromised,
particularly for a discretization without additional stabilization
terms.

\begin{table}
  \centering
  \caption{Nonlinear and linear solver performance for MHD generator
    example.  Problems with 724K DoFs were run on 5 cores, those with
    5.6M DoFs were run on 40 cores (1 node), and those with 44M DoFs
    were run on 400 cores (10 nodes).}
  \label{tab:generator}
  \begin{tabular}{c||cc|ccc|cc}
    \toprule
    Finest grid size & \multicolumn{2}{c|}{724,725}
    & \multicolumn{3}{c|}{5,600,229}
    & \multicolumn{2}{c}{44,023,749}\\
    Number of levels & 2 & 3 & 2 & 3 & 4 & 3 & 4 \\
    \midrule
    Newton Steps & 5 & 5 & 4 & 4 & 5 & 4 & 4 \\
    Total Linear iterations & 52 & 62 & 42 & 44 & 117 & 58 & 184 \\
    Solution Time (minutes) & 4.12 & 4.31 & 5.09 & 5.18 & 9.86 & 5.68
    & 12.85 \\
    \bottomrule
    \end{tabular}
  \end{table}

For the smallest problems, memory usage peaked at approximately 4 GB
per core (20GB total), which rose slightly for the middle problem sizes, which
required approximately 170 GB total, filling the available memory on
one node of Niagara (after accounting for required system overhead).  For the finest problem size, there was
insufficient memory per node to run the solvers on 8 nodes, so
these results were run on 10 nodes, requiring about 160 GB/node for
the 4-level solver (with increased memory costs due to GMRES storage
for the large iteration counts) and 153 GB/node for the 3-level
solver.  Solution times scale relatively clearly with linear iteration
counts, leading to roughly scalable 3-level solvers (although we note
the finest grid results clearly benefit in this comparison from the
added CPUs available on 10 nodes).  Finally, we note that it was not possible to
test a ``deeper''
multigrid hierarchy on the finest mesh,
due to a software requirement that there be at least
one element per MPI thread on the coarsest mesh of the hierarchy.  The
coarsest mesh for the smaller problem sizes contains 240 elements (a
$10 \times 2 \times 2$ hexahedral mesh refined into tetrahedra), which
does not satisfy this constraint.

\section{Conclusions}

We present a monolithic multigrid preconditioner for the
linearization of a mixed finite element discretization of the
equations of viscoresistive incompressible MHD.  A distinctive
feature of this
formulation is the presence of two Lagrange multipliers, one to
enforce the incompressibility of fluid velocity, and one to enforce a
similar divergence-free constraint on the magnetic field.  Three
variants of Vanka relaxation are considered; while the segregated
approach breaks down when coupling between the fluid velocity and
magnetic field is significant, both the purist and coupled approaches
appear to provide robust solvers, with the coupled approach offering the
most efficient performance.  The associated multigrid solver is
effective for the stationary and transient cases up to 170 million
degrees of freedom.

The clear limitation of the approach presented here is the low order
and lack of stabilization for the finite element
discretization considered, and extending this to a better
finite element discretization is the main focus for future work.  We
note the preconditioners proposed in \cite{EGPhillips_etal_2016a,
  shadid2016scalable} are applied to stabilized finite element
discretizations similar to those considered here, while excellent
augmented Lagrangian preconditioners are developed for stabilized discretizations 
of the Navier--Stokes at
large Reynolds numbers in
\cite{farrell2018augmented, farrell2020}.  Another clear avenue for
future work is the development of similar monolithic multigrid
approaches for more detailed models of plasma, such as those considered
in \cite{srinivasan2011analytical,loverich2011discontinuous,meier2012general,miller2019imex,phillips2020enabling}.

\bibliographystyle{siam}
\bibliography{MGrefs}

\end{document}